\documentclass[reqno,10pt]{amsart}
\usepackage{amsthm}
\usepackage{amsmath}
\usepackage{amssymb}
\usepackage{amsfonts}
\usepackage{mathrsfs}
\usepackage[pdftex]{graphicx}

\usepackage[utf8]{inputenc}
\usepackage{latexsym}
\usepackage{verbatim}

\usepackage[all]{xy}
\usepackage[bookmarksnumbered,colorlinks,urlcolor=blue]{hyperref}
\usepackage{graphicx}
\usepackage{graphics}
\usepackage{float}
\usepackage{enumerate}
\def\bb#1\eb{\textcolor{blue}
{#1}} %
\def\br#1\er{\textcolor{red}
{#1}} %
\newcommand{\qcd}{\begin{flushright} $\Box$ \end{flushright}}

   \def\br#1\er{\textcolor{red}{#1}} %
      \def\bb#1\eb{\textcolor{blue}{#1}} %

\title[]{A note on causality conditions on covering spacetimes}

\author[E. Minguzzi]{Ettore Minguzzi}
\address{Dipartimento di Matematica e Informatica ``U. Dini'',
\hfill\break\indent Universit\`{a} degli Studi di Firenze,
\hfill\break\indent Via S. Marta 3, I-50139,
\hfill\break\indent Firenze, Italy}
\email{ettore.minguzzi@unifi.it}

\author[I.P. Costa e Silva]{Ivan P. Costa e Silva}
\address{Department of Mathematics, Universidade Federal de Santa Catarina,
	\hfill\break\indent 88.040-900
	Florianópolis-SC, Brazil.}
\email{pontual.ivan@ufsc.br}


\begin{document}
\newtheorem{thm}{Theorem}[section]
\newtheorem{prop}[thm]{Proposition}
\newtheorem{lemma}[thm]{Lemma}
\newtheorem{cor}[thm]{Corollary}
\theoremstyle{definition}
\newtheorem{defi}[thm]{Definition}
\newtheorem{notation}[thm]{Notation}
\newtheorem{exe}[thm]{Example}
\newtheorem{conj}[thm]{Conjecture}
\newtheorem{prob}[thm]{Problem}
\newtheorem{rem}[thm]{Remark}
\newtheorem{conv}[thm]{Convention}
\newtheorem{crit}[thm]{Criterion}
\newtheorem{claim}[thm]{Claim}

\newcommand{\ben}{\begin{enumerate}}
\newcommand{\een}{\end{enumerate}}

\newcommand{\bit}{\begin{itemize}}
\newcommand{\eit}{\end{itemize}}

\begin{abstract}
A number of techniques in Lorentzian geometry, such as those used in the proofs of singularity theorems, depend on certain smooth coverings retaining interesting global geometric properties, including causal ones.
In this note we give explicit examples showing that, unlike some of the more commonly adopted rungs of the causal ladder such as strong causality or global hyperbolicity, less-utilized conditions such as causal continuity or causal simplicity {\em do not} in general pass to coverings, as already speculated by one of the authors (EM). As a consequence, any result which relies on these causality conditions transferring to coverings must be revised accordingly. In particular, some amendments in the statement and proof of a version of the Gannon-Lee singularity theorem previously given by one of us (IPCS) are also presented here that address a gap in its original proof, simultaneously expanding its scope to spacetimes with lower causality.
\end{abstract}

\maketitle

\vspace*{-.5cm}

\section{Introduction}\label{section1}

Let $(M^n,g)$ be a {\em spacetime}, i.e., a pair consisting of a connected $C^{\infty}$ smooth manifold (Hausdorff and secound-countable) $M$ of dimension $n\geq 2$ and a time-oriented $C^\infty$ Lorentz metric $g$. Let $\pi:\tilde{M}\rightarrow M$ be a smooth covering map, and endow $\tilde{M}$ with the pullback metric $\tilde{g}:= \pi ^{\ast}g$ and the induced time-orientation. Since $\pi$
is a local isometry, any {\em local} geometric condition that might hold on $(M,g)$, such as (say) $Ric_g(v,v)\geq 0$ for all lightlike vectors $v \in TM$, must hold on $(\tilde{M},\tilde{g})$ as well.

Geodesic (lightlike, or timelike, or spacelike) completeness, on the other hand, is a key example of a {\em global} geometric feature which holds on $(M,g)$ if and only if it holds on $(\tilde{M},\tilde{g})$. This is of substantial technical importance in the proofs of singularity theorems, since some constructions are carried out on (a suitable choice of) $(\tilde{M},\tilde{g})$ \cite{C,G,Haag,H,HE}. For example, (cf. \cite[Prop.\ 14.2]{oneill}) if $(M,g)$ admits a (topologically) closed connected spacelike hypersurface $\Sigma \subset M$, then $(\tilde{M},\tilde{g})$ can be chosen so that it possesses a diffeomorphic copy $\tilde{\Sigma}$ of $\Sigma$ which is in addition {\em acausal},
so that its Cauchy development can be suitably analyzed for the existence of certain maximal geodesics normal to $\tilde{\Sigma}$, an important step in some proofs. Of course, such constructions can only be meaningfully carried out provided there is some
control on whether the required properties still hold on covering manifolds.

An important hypothesis in a number of theorems, and in singularity theorems in particular, is on which rung of the so-called {\em causal ladder} \cite{minguzzi_living_reviews,causalladder} of spacetimes $(M,g)$ sits. In view of the remarks in the previous paragraph, it is therefore of interest to know whether that rung is shared with $(\tilde{M},\tilde{g})$. A positive statement in this regard has been summarized by one of us (EM) as follows. (See \cite[Thm.\ 2.99]{minguzzi_living_reviews} for an extended discussion and a proof.)

\begin{thm}\label{thm1}
Let $\pi : (\tilde{M},\tilde{g})\rightarrow (M,g)$ be a Lorentzian covering. If $(M,g)$ is chronological [resp. causal, non-totally
imprisoning, future/past distinguishing, strongly causal, stably causal, globally
hyperbolic] then $(\tilde{M},\tilde{g})$ has the same property.
\end{thm}

\noindent Just after the statement of that theorem, the author mentions in passing that {\em reflectivity} and {\em closure of causal futures/pasts} of points in $(M,g)$ do not seem to pass to coverings. It is the purpose of this note to both confirm the latter claim by means of concrete (counter)examples, as well as to modify the statement and proof of \cite[Thm.\ 2.1]{costa_e_silva_gannon} to incorporate this discovery. As it stands, the latter proof has a gap if the underlying spacetime is not globally hyperbolic, precisely because it assumes without further discussion that simple causality also applies to a certain Lorentzian covering thereof.\footnote{IPCS wishes to thank Roland Steinbauer and Benedict Schinnerl for calling his attention to the fact that that assumption is made in \cite{costa_e_silva_gannon} without proper justification.}

The version we present here, however, is {\em not} a mere amendment.
Recently it has been shown that the assumptions of Penrose's theorem can be improved by weakening global hyperbolicity to past reflectivity \cite{minguzzi_new}. Adapting some arguments to the Gannon-Lee case we are able to accomplish a similarly interesting result, that is,
we can dramatically decrease the causality requirements of the Gannon-Lee theorem by demanding $(M,g)$  {\em and its coverings} to be just past reflecting. Although the latter criterion might seem impractical, the (proof of) \cite[Thm.\ 4.10]{minguzzi_living_reviews} actually implies that {\em any} spacetime - irrespective of added causal assumptions - with a past-complete conformal timelike Killing vector field satisfies it.

This paper is organized as follows. In section \ref{s1}, we briefly discuss the notions of past/future reflectivity, an in particular the result mentioned towards the end of the previous paragraph. Then we present the two central examples, one of which is inspired by
the spacetime constructed by Hedicke and Suhr in \cite[Thm.\ 2.7]{hedicke19} (which was constructed with the entirely different purpose of providing an example of a causally simple spacetime for which the space of null geodesics is not Hausdorff). These examples show that $(i)$ the closure of the causal relation, and $(ii)$ past reflectivity, do not pass to $(\tilde{M},\tilde{g})$ while holding on $(M,g)$.
In section \ref{s2} we give an alternative statement and proof of the Gannon-Lee theorem presented in \cite[Thm.\ 2.1]{costa_e_silva_gannon} in terms of past reflectivity of Lorentzian coverings.

We shall assume throughout that the reader is familiar with the elements of causal theory in the core references \cite{beem,oneill}, up to and including the best-known singularity theorems originally proven by R. Penrose and S.W. Hawking described in those references. We also assume the reader is acquainted with the basic structure of the causal ladder, that is, the basic hierarchy of causal conditions listed from the weakest - {\em non-totally viciousness} - at the bottom, through the strongest - {\em global hyperbolicity} at the top, which can be found, e.g., in Ch. 2 of \cite{beem}. However, as mentioned above, since the notions of past/future reflectivity are somewhat less known, we recall these in section \ref{s1}.

The results in this paper are purely geometric in that no field equations are assumed, and they hold for any spacetime dimension $\geq 3$. Our conventions for the signs of spacetime curvature are those of \cite{beem}, but for the mean curvature vector of submanifolds within it we use those of \cite{oneill}. In particular, we borrow the definition of {\em convergence} of a semi-Riemannian submanifold from the latter reference (cf. \cite[Def.\  10.36]{oneill}). Unless otherwise explicitly stated, all maps and (sub)manifols are assumed to be $C^{\infty}$, and submanifolds to be embedded.

\section{(Counter)examples of causality conditions on coverings}\label{s1}

As announced in the Introduction, in this section we discuss some facts about reflectivity and present examples showing that certain causal conditions which hold on a spacetime $(M,g)$ may fail to hold on one (or more) of its coverings, including the universal covering.

First, we recall some terminology.
\begin{defi} \label{reflectdef}
A spacetime $(M,g)$ is said to be {\em past reflecting} [resp. {\em future reflecting}] if one (and hence both) of the two equivalent statements hold for any two $p,q\in M$:
\begin{itemize}
\item[i)] $I^+(q)\subset I^+(p) \Rightarrow I^-(p)\subset I^-(q)$ [resp. $I^-(q)\subset I^-(p) \Rightarrow I^+(p)\subset I^+(q)$];
\item[ii)] $q \in \overline{I^+(p)} \Rightarrow p \in \overline{I^-(q)} $ [resp. $q \in \overline{I^-(p)} \Rightarrow p \in \overline{I^+(q)} $].
\end{itemize}
If $(M,g)$ is {\em both} past and future reflecting, it is simply said to be {\em reflecting}.
\end{defi}
\noindent For a more detailed discussion on reflectivity, and a proof of the equivalence of $(i)$ and $(ii)$ in Def. \ref{reflectdef} (as well as other equivalent statements) see, e.g., \cite[Section\ 4.1]{minguzzi_living_reviews}.

Unlike usual causality conditions, reflectivity {\em by itself} lies outside the causal ladder, and is known as a {\em transversal} condition. (See, for example, the discussion around \cite[Fig.\ 2]{k_causalityminguzzi} or \cite{minguzzi_living_reviews}.)

There is a well-known relationship between past [resp. future] reflectivity and a so-called past [resp, future ] {\em volume function} $t^{-}$ [resp. $t^{+}$] (cf. \cite[Prop.\ 3.21]{beem}). From our perspective here, however, there is another useful sufficient condition ensuring that not only a spacetime $(M,g)$ {\em but also its coverings} are (past/future) reflecting. Recall that a smooth vector field $X:M\rightarrow TM$ is a conformal Killing vector field if $\mathcal{L}_Xg = \sigma \cdot g$, where $\mathcal{L}$ denotes the Lie derivative and $\sigma \in C^{\infty}(M)$. (Of course, $X$ is Killing if and only if $\sigma \equiv0$.)
\begin{prop}\label{criterion}
If a spacetime $(M,g)$ admits a past-[resp. future-]complete conformal timelike Killing vector field, then $(M,g)$   is past [resp. future ] reflecting. In addition, given any Lorentz covering $\pi : (\tilde{M},\tilde{g})\rightarrow (M,g)$, the spacetime (with the induced time-orientation) $(\tilde{M}, \tilde{g})$ is also past [resp. future] reflecting.
\end{prop}
\noindent {\em Proof.} We just deal with the past complete case, since the future case follows by time-duality. Let $X\in \mathfrak{X}(M)$ be a past-complete conformal timelike Killing vector field, and write $\beta:=\sqrt{|g(X,X)|}>0$. It is well-known, and in any case easy to check, that $X$ is a {\em Killing} vector field for the spacetime $(M,\hat{g})$, where $\hat{g}:= \beta ^{-2}g$. But past reflectivity, just as any other causal property, is conformally invariant, and hence holds on $(M,g)$ if and only if it holds on $(M,\hat{g})$. Therefore, there is no loss of generality in assuming that $X$ is in fact a Killing vector field. The proof of \cite[Thm.\ 4.10]{minguzzi_living_reviews} then establishes that $(M,g)$ is past reflecting. (In fact, this proof works directly just as well for a conformal Killing vector field, providing an alternative to the argument above.)  In order to complete the proof, just note that given any Lorentz covering $\pi : (\tilde{M},\tilde{g})\rightarrow (M,g)$, by standard smooth covering properties there exists a unique vector field $\tilde{X} \in \mathfrak{X}(\tilde{M})$ such that
\begin{equation}\label{adhoc}
d\pi _{\tilde{p}}(\tilde{X}_{\tilde{p}})=X_{\pi(\tilde{p})}, \quad \forall \tilde{p}\in \tilde{M}.
\end{equation}
Since $\pi$ is a local isometry, $\tilde{X}$ is a
(conformal) Killing vector field if and only if $X$ is also. Finally, (\ref{adhoc}) means that any integral curve of $\tilde{X}$ is a lift through $\pi$ of an integral curve of $X$, so $\tilde{X}$ is also past-complete.
\qcd

Proposition \ref{criterion} shows, in particular, that reflectivity may well be present even in non-chronological spacetimes. For example, any {\em compact} spacetime (which can never be chronological - cf. \cite[Lemma\ 14.10]{oneill}) endowed with a (necessarily complete) conformal Killing vector field is reflecting. As a famous example, G\"odel spacetime admits a complete timelike Killing vector field (see, e.g. \cite[Section\ 5.7, p.\ 168]{HE} for a brief discussion), and thus it is also reflecting. (Since it is diffeomorphic to $\mathbb{R}^4$ it is its own universal covering.) Another simple example is given by identifying $n$-dimensional Minkowski spacetime with, say, global coordinates $(t,x^1,\ldots,x^{n-1})=:(t,x)$,  along $t$, e.g., identifying $(t,x)\sim (t+1,x)$. However, if combined with other, even quite mild, causality conditions, reflectivity does imply - and is implied by - some rather strong causality requirements \cite{minguzzi_living_reviews}:
\begin{itemize}
\item[a)] $(M,g)$ is distinguishing $+$ reflecting $\Leftrightarrow$ $(M,g)$ is causally continuous;
\item[b)] $(M,g)$ is causal $+$ closure of causal relation $\Leftrightarrow$ $(M,g)$ is causally simple;
\item[c)] closure of causal relation $\Rightarrow$ $(M,g)$ is reflecting;
\item[d)] $(M,g)$ is non-totally vicious + reflecting $\Rightarrow$ $(M,g)$ is chronological.
\end{itemize}
(In particular, item $(d)$ implies, in view of Prop. \ref{criterion}, that any non-chronological spacetime - for instance a compact spacetime - with a complete conformal timelike Killing vector field is necessarily totally vicious.)

We are ready to discuss our main examples.

\begin{exe}\label{exe1}
We construct a causally simple $3$-dimensional spacetime $(M^3,g)$ for
which the {\em universal covering} $(\tilde{M}, \tilde{g})$ is such that $\tilde J$,
the causal relation on $\tilde M$, is not closed. Thus, $(\tilde{M}, \tilde{g})$ cannot be causally simple.

Consider the static spacetime
\[
(M,g)=(\mathbb{R}\times \Sigma, -d t^2+ \kappa),
\]
where $\partial _t$ is taken to be future-directed and $(\Sigma,\kappa)$ is some Riemannian manifold.  (Note that since $\partial _t$ is a complete Killing vector field, such a spacetime is in particular reflecting by Prop. \ref{criterion}.) Hedicke and Suhr have shown (cf. remark just after \cite[Cor.\ 4.2]{hedicke19}) that $(M,g)$ is causally simple if and only if $(\Sigma,\kappa)$ is {\em geodesically convex}, i.e. any two $x,y\in \Sigma$ are connected in $\Sigma$ by a (not necessarily unique) {\em minimizing} $\kappa$-geodesic.

Of course, if $(\Sigma,\kappa)$ is a {\em complete} Riemannian manifold, then by (the proof of) Hopf-Rinow theorem it is necessarily geodesically convex. (Indeed it is well-know that $(M,g)$ is actually globally hyperbolic in this case - cf., e.g., \cite[Thm.\ 3.67]{beem}.) The key to our example is therefore to find an incomplete Riemannian manifold which is geodesically connected in this sense, but that has some Riemannian covering that isn't.

Concretely, we choose $(\Sigma,\kappa)$ to be a suitable surface of revolution in $\mathbb{R}^3$  given as follows. Consider the map
\[
\varphi : (0,+\infty)\times \mathbb{R} \to \mathbb{R}^3, \quad (u,v) \mapsto
(A u\cos\ v, A u\sin \ v, u),
\]
where $A$ is for now some positive number which we shall further restrict momentarily. $(\Sigma, \kappa)$ is then the image of $\varphi$ endowed with the (Riemanniann) metric induced by the Euclidean metric on $\mathbb{R}^3$. This image is of course a (flat) cone without its apex, and its axis lies along the $z$-axis.

The cone $\Sigma$ can be isometrically identified with a circular sector $C$ on the plane $\mathbb{R}^2$ whose sides are identified. This fact greatly simplifies the analysis (cf. Figure \ref{sec}).

The angle comprehended by this sector is
\[
\theta:=\frac{2A\pi}{\sqrt{1+A^2}}.
\]
Provided $\theta<\pi$ (equivalently $A<1/\sqrt{3}$) any two points of $C$ are connected by a minimal geodesic which is a segment not intersecting the singularity. As a consequence $(\Sigma, \kappa)$ is geodesically convex. (Roughly, no two points $x,y \in \Sigma$ can be connected by a minimal geodesic passing through the singularity, for the sum of the Euclidean lengths of the segments connecting $x,y$ to the origin is always larger than the segment in $C$ connecting them.) Thus, the latter bound on $A$ guarantees the causal simplicity property for $(M,g)$ by the Hedicke-Suhr criterion.

Choose now two points $a,b\in C$ as follows. $a$ is chosen on one of the sides of $C$, and $b$ is chosen  on the bisecting line of $\theta$. Consider the fixed-endpoint homotopy class of a curve $\gamma_1$ starting from $a$ and reaching $b$, that revolves once in the positive direction over the singularity at $x=0$. Denote this class by $[\gamma_1]$, which is obviously distinct from that of any minimal geodesic connecting them (there are two minimal geodesics in total, one being  denoted $\sigma$ in  Figure \ref{sec}).
In Figure \ref{sec} we draw two adjacent copies of $C$; it is then easily seen that provided $\theta>2\pi/3$ (equivalently $A>\frac{1}{2\sqrt{2}}$, so overall $A=1/2$ nd $\theta\simeq 161^\circ$ is a possible choice),
we have
\[
\ell _0:= \textrm{inf}_{\gamma\in [\gamma_1]}
\ell^\kappa(\gamma)=\ell^\kappa(\eta_1)+\ell^\kappa(\eta_2)>
\ell^\kappa(\sigma)
\]
where $\ell^\kappa$ is the length functional in $(\Sigma,\kappa)$ and
where $\eta_1$ and  $\eta_2$ are incomplete geodesics that connect $a$
to the singularity and the singularity to $b$, respectively. The last inequality follows from the triangle inequality of the Euclidean plane.
 Notice that the infimum of length $\ell _0$ is not realized since the putative `curve'
$\eta_2 \circ \eta_1$ would pass through the singularity if connected.

\begin{figure}[ht]
\centering
\includegraphics[width=8cm]{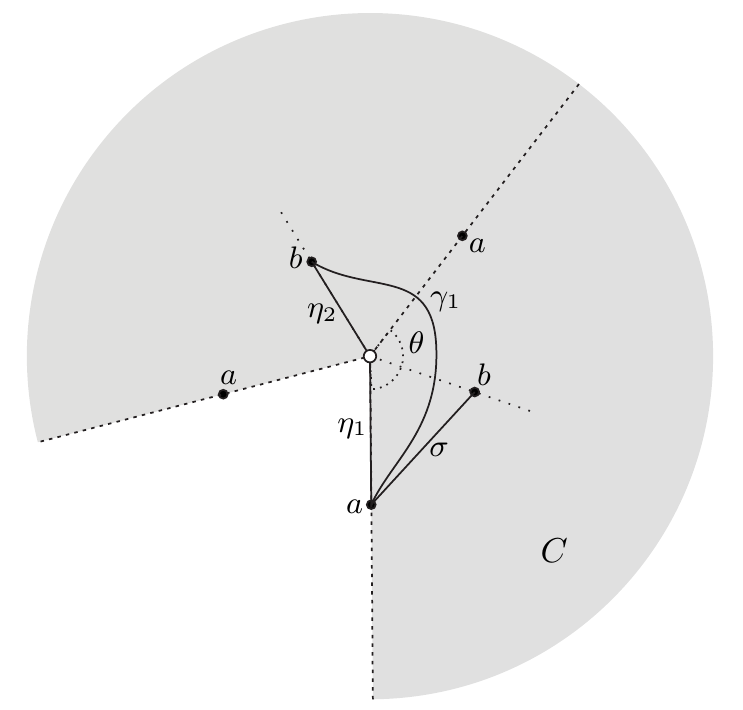}
\caption{Two adjacent copies of the sector $C$ on the Euclidean plane, where we unwrap the curve $\gamma_1$. By identifying suitable dotted lines we recover our cone $\Sigma$. Thinking of $\gamma_1$ as an elastic band it is clear than its shorter length would be assumed in the configuration $\eta_2\circ \eta_1$ provided $3\theta/2>\pi$. Similarly, the minimal curve connecting $a$ to itself does not intersect the singularity provided $\theta<\pi$. The full depicted region can be thought as a portion of the covering $\tilde \Sigma$.
} \label{sec}
\end{figure}

Now, the universal cover $(\tilde{M},\tilde{g})$ of $(M,g)$ is clearly $(\mathbb{R}\times \tilde{\Sigma},-dt^2+ \tilde{\kappa})$, where $(\tilde{\Sigma}, \tilde{\kappa})$ is the universal (Riemannian) covering of $(\Sigma, \kappa)$. The construction above actually shows that a given copy of $a$ and other copies of $b$ on $\tilde{\Sigma}$ not projecting onto the same sector are not connected by a minimizing geodesic. We conclude that $(\tilde{\Sigma},\tilde{\kappa})$ is not geodesically connected, so the causal relation $\tilde{J}$ on $(\tilde{M},\tilde{g})$ is not closed by the Hedicke-Suhr criterion.
\end{exe}

We now present an example of causally continuous (and in particular reflecting) spacetime $(M,g)$ for which the universal covering $(\tilde M, \tilde g)$ is not past reflecting.
\begin{exe}\label{exe2}
Start with Minkowski 2d spacetime $\mathbb{R}^2$ of
coordinates $(t,x)$, endowed with the metric $-d t^2+d x^2$. Let
$r=(0,0)$ and $r_k=(0,-1/k)$ ($k\in \mathbb{N}$). Consider now the spacetime given as the set
$M=\mathbb{R}^2 \backslash \{r,r_1,r_2, \cdots \}$ endowed with the restricted
metric and time-orientation, see Figure \ref{ref}. This spacetime is causally continuous, as can be easily established by
looking at the continuity of the volume functions \cite[Def.\
4.6(vii)]{minguzzi_living_reviews}.  Let $p=(-1,1)$ and $q=(1,-1)$, so that $q\in
\overline{I^+(p)}$ and $p\in \overline{I^-(q)}$.

Consider the universal covering $\pi:(\tilde M,\tilde{g}) \rightarrow (M,g)$. Let $\gamma:[0,1]\rightarrow M$ be any curve such that $\gamma(0)=p$, $\gamma(1)=q$, and $\gamma (0,1) \subset I^+(p)$, and consider its lift $\tilde{\gamma}$ starting at a fixed representative $\tilde p \in \pi ^{-1}(p)$. Then, for each $k\in \mathbb{N}$ with $k\geq 2$,
\[
\tilde{q}_k:=\tilde{\gamma}(1-1/k) \in \pi^{-1}(I^+(p))\equiv \tilde{I}^+(\pi^{-1}(p));
\]
but since the restrictions $\gamma |_{[0,1-1/k]}$ of the curve $\gamma$ to the intervals $[0,1-1/k]$ are clearly endpoint homotopic to any future-directed timelike curve from $p$ to $\gamma (1-1/k)$ we conclude that $\tilde{q}_k \in \tilde{I}^+(\tilde{p})$, so $\tilde{q}:= \tilde{\gamma}(1) \in \pi ^{-1}(q)\cap \overline{\tilde I^+(\tilde p)}$.

\begin{figure}[ht]
\centering
\includegraphics[width=9cm]{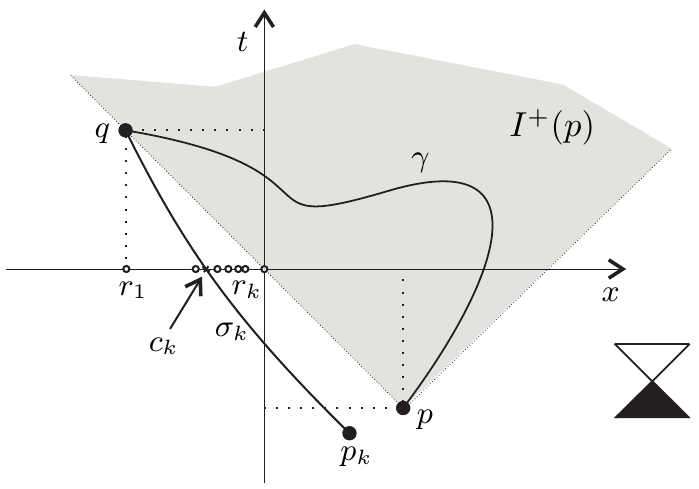}
\caption{The construction of Example \ref{exe2} on $(M,g)$. Note that the $r_k$'s are not points of the spacetime.
} \label{ref}
\end{figure}

We claim, however, that $\tilde p \notin \overline{\tilde I^-(\tilde q)}$, i.e., past reflectivity is violated on $(\tilde M,\tilde{g})$. For suppose there is a sequence $(\tilde{p}_k)\subset \tilde {I}^-(\tilde q)$ converging to $\tilde{p}$. Pick any future-directed timelike curves $\tilde{\sigma}_k:[0,1]\rightarrow \tilde{M}$ starting at $\tilde{p}_k$ and ending at $\tilde{q}$, so that $p_k:=\pi(\tilde{p}_k) \rightarrow p$. The curves $\sigma_k:= \pi\circ \tilde{\sigma}_k$ are timelike, so $(p_k)\subset I^-(q)$. Moreover, the intersections of the images of the $\sigma_k$'s with the axis $t=0$ clearly must occur at points $c_k \rightarrow (0,0)$ (in $\mathbb{R}^2$). Therefore, for some subsequence $(c_{k_i})_{i\in \mathbb{N}}$ we can assume that each $c_{k_i}$ belongs to a different segment $\{0\}\times (-\frac{1}{k_i}, -\frac{1}{k_i+1})$. Pick any distinguished small open disc $U\ni p$ in $M$ not intersecting the $t=0$ axis. We can assume, without loss of generality, that $p_{k_i} \in U$ for every $i\in \mathbb{N}$. But then the points $\tilde p_{k_i}$ must be in different connected components of $\pi^{-1}(U)$, and  thus we cannot have $\tilde p_{k_i}\rightarrow \tilde{p}$, a contradiction.
\end{exe}

We mention in passing that the properties  `reflectivity' and  `closure of the causal relation' not only do not pass to coverings, but are also known to be distinguished among all causality properties for not being preserved by {\em isocausal mappings}. We do not know whether this has any deeper significance or is just a coincidence. 
We do not pursue this point any further. The interested reader is referred to \cite{garciaparrado05b} for definitions and details.

\section{A new version of the Gannon-Lee theorem in low causality} \label{s2}

The {\em Gannon-Lee singularity theorem} was independently discovered by D. Gannon and C.W. Lee in 1975/1976 \cite{gannon1,gannon2,lee}. Its importance lies in its application to general relativity, wherein it suggests that certain ``localized non-trivial topological structures'' in spacetime (meaning non-trivial fundamental group of certain spacelike hypersurfaces), such as wormholes, are {\em gravitationally unstable}, at least if one neglects quantum effects.

To make precise statements, we again fix an $n$-dimensional spacetime $(M,g)$, but throughout this section we assume $n\geq 3$. We shall recall some terminology which we believe is unfamiliar for many readers, largely following \cite{costa_e_silva_gannon}.

Fix a smooth, connected, spacelike partial Cauchy hypersurface\footnote{Recall that a {\em partial Cauchy hypersurface} is by definition an acausal edgeless subset of a spacetime, which means in particular that it is a topological (i.e. $C^0$) hypersurface \cite{oneill}. In this paper, however, we always deal with smooth hypersurfaces.} (i.e, a submanifold of codimension one) $\Sigma^{n-1} \subset M$, and a smooth, connected, compact spacelike submanifold $S^{n-2} \subset M$ of codimension two.

 Suppose $S$ {\em separates} $\Sigma$, i.e., $S \subset \Sigma$ and $\Sigma \setminus S$ is not connected. This means, in particular, that $\Sigma \setminus S$ is a disjoint union $\Sigma_{+}\dot{\cup} \Sigma_{-}$ of open submanifolds of $\Sigma$ having $S$ as a common boundary. We shall loosely call $\Sigma_{+}$ [resp. $\Sigma_{-}$] the {\em outside} [resp. {\em inside}] of $S$ in $\Sigma$. (In most interesting examples there is a natural choice for these.) It also means that there are unique unit spacelike vector fields $N_{\pm}$ on $S$ normal to $S$ in $\Sigma$, such that $N_{+}$ [resp. $N_{-}$] is outward-pointing [resp. inward-pointing].

Let $U$ be the unique timelike, future-directed, unit normal vector field on $\Sigma$. Then $K_{\pm} := U|_{S}+N_{\pm}$ are future-directed null vector fields on $S$ normal to $S$ in $M$ spanning the normal bundle $NS \subset TM$ (which is in particular trivial). The {\em outward }[resp. {\em inward}] {\em null convergence} of $S$ in $M$ is the smooth function $k_{+}: S \rightarrow \mathbb{R}$ [resp. $k_{-}: S \rightarrow \mathbb{R}$] given by
\begin{equation}
\label{maineq1}
k_{+}(p) = \langle H_p, K_{+}(p)\rangle _{p} \mbox{  [resp. $k_{-}(p) = \langle H_p, K_{-}(p) \rangle_p $]},
\end{equation}
for each $p \in S$, where $H_p$ denotes the mean curvature vector of $S$ in $M$ at $p$ \cite{oneill}, and we denote $g$ as $\langle \, , \, \rangle$ here and hereafter, if there is no risk of confusion. Under the sign conventions we adopt here, if $S$ is a round sphere in a Euclidean slice of Minkowski spacetime with the obvious choices of inside and outside, then $k_{+}<0$ and $k_{-}>0$. One also expects this to be the case if $S$ is a ``large"  sphere in an asymptotically flat spatial slice.

Using the terminology above, we shall adopt the following useful definition:
\begin{defi}[Asymptotically regular hypersurface]
 \label{asympticallyregular}
 A smooth, connected, spacelike partial Cauchy hypersurface $\Sigma \subset M$ is {\em asymptotically regular} if there exists a smooth, connected, compact submanifold $S \subset \Sigma$ of dimension $n-2$ such that
\begin{itemize}
\item[i)] $S$ separates $\Sigma$, and $\overline{{\Sigma}}_{+} \equiv S \cup \Sigma_{+}$ is non-compact;
\item[ii)] The map $h_{\#}: \pi_1(S) \rightarrow \pi_1(\overline{\Sigma}_{+})$ induced by the inclusion $h: S \hookrightarrow \overline{\Sigma}_{+}$ is onto;
\item[iii)] $S$ is {\em inner trapped}, i.e., $k_{-} >0$ everywhere on $S$.
 \end{itemize}
we shall call such an $S$ an {\em enclosing surface} in $\Sigma$.
\end{defi}
Let us briefly pause to explain the motivation behind the clauses $(i)-(iii)$ of this definition. First, it is meant as a convenient adaptation of Gannon's definition of a {\em regular near infinity} hypersurface, so item $(i)$ presents no novelty. Clause $(ii)$, however, might look somewhat opaque.  But it simply means that the (closure of the) outside of $S$ has only topological (or more precisely path-homotopic) complexities arising from having $S$ itself as a boundary. Specifically, {\em it means that every loop in the exterior of $S$ in $\Sigma$ is homotopic to a loop on $S$}. Note that this is certainly the case if $\overline{\Sigma}_{+} \equiv S \cup \Sigma_{+}$ is homeomorphic to $S \times [0, +\infty)$, as in the original Gannon-Lee theorem, but the condition as stated gives rise to the much wider set of topological possibilities which are likely to arise in higher dimensions. Note that $(iii)$ refers only to the inward-pointing family of null geodesics normal to $S$, namely that they converge ``on average''.

It should be clear, furthermore that while conditions $(i)-(iii)$ are naturally expected to occur in asymptotically flat spatial slices of a given spacetime, they are actually much weaker requirements than asymptotic flatness, and can occur for most of the falloff conditions for initial data slices studied in the extant literature.

Another important technical notion introduced in \cite{costa_e_silva_gannon} is that of a {\em piercing} for a partial Cauchy hypersurface $\Sigma$.
 \begin{defi}[Piercing] \label{piercing}
 We say that a smooth future-directed timelike vector field $X:M \rightarrow TM$ is a {\em piercing} of the partial Cauchy hypersurface $\Sigma$ (or {\em pierces} $\Sigma$) if every maximally extended integral curve of $X$ intersects $\Sigma$ exactly once.
 \end{defi}
If $X\in \mathfrak{X}(M)$ is a piercing for $\Sigma$, there is no loss of generality in assuming that $X$ is a complete vector field, and we shall do so in what follows. Using its flow, it is not difficult to show that if such a piercing exists, then $M$ is diffeomorphic to $\mathbb{R}\times \Sigma$.

 Of course, an asymptotically regular partial Cauchy hypersurface $\Sigma$ and/or a piercing of $\Sigma$ may not exist for general spacetimes. On the other hand, if $(M,g)$ is globally hyperbolic and $\Sigma$ is a Cauchy hypersurface, then {\em every} smooth future-directed timelike vector field in $M$ pierces $\Sigma$. {\em However, the existence of a piercing for a partial Cauchy hypersurface $\Sigma$ is strictly weaker than the requirement that $\Sigma$ be Cauchy}. In \cite{costa_e_silva_gannon}, anti-de Sitter spacetime is given as an example, but since the latter is causally simple, it is instructive to give another class of simple extra examples of (slightly) lower causality in preparation for our main theorem.

\begin{exe}\label{morexe1}
Start with $n$-dimensional Minkowski spacetime
\[
(\mathbb{R}^n, -dt^2 +\sum _{i=1}^{n-1}(dx^i)^2)
\]
with $n\geq3$ and the usual time-orientation, and denote $z=(x^1,\ldots,x^{n-1})\in \mathbb{R}^{n-1}$, and a generic spacetime point by $(t,z)$. Fix $k$ pairwise distinct points $z_1,\ldots,z_k \in \mathbb{R}^{n-1}$ ($k\geq 1$), and pick any $R>0$ for which the sphere $\mathbb{S}^{n-2}(R) \subset \mathbb{R}^{n-1}$ of radius $R$ and center at the origin is such that $z_1,\ldots,z_k$ all lie in the interior of that sphere inside $\mathbb{R}^{n-1}$. Finally, consider the $k$ timelike lines
\[
\ell_{j}:= \{(t,z_j)\, : \, t\in \mathbb{R}\}, \quad j=1,\ldots,k,
\]
and define $(M^n,g)$ as the flat spacetime given by $M=\mathbb{R}^n \setminus (\cup _{j=1}^k \ell_{j})$   with $g$ the restricted metric and time orientation. The restriction $X$ of $\partial _t$ to $M$ is a complete Killing vector field for $(M,g)$, so by Prop.\ \ref{criterion}  not only $(M,g)$ but any of its Lorentzian coverings is reflecting. Since $(M,g)$ is of the form $(\mathbb{R}\times \hat{\Sigma},-dt^2+\kappa)$ as in Example \ref{exe1} with $\hat{\Sigma} := \mathbb{R}^{n-1}\setminus \{z_1,\ldots,z_k\}$ and $\kappa$ the flat metric thereon, $(\hat{\Sigma},\kappa)$ is clearly not geodesically convex, so $(M,g)$ is {\em not} causally simple by the Hedicke-Suhr criterion, although it is causally continuous.

The partial Cauchy hypersurface $\Sigma:= \{0\}\times \hat{\Sigma}$ is asymptotically regular with enclosing surface $S:= \{0\}\times \mathbb{S}^{n-2}(R)$. Note, furthermore, that $\Sigma$ is pierced by $X$. However, when $k\geq 2$ and $n=3$ there are curves in $\Sigma$ with endpoints on $S$ which cannot be deformed therein to a curve on $S$, so the inclusion-induced homomorphism $\pi_1(S)\hookrightarrow \pi_1(\Sigma)$ is {\em not} surjective in this case. This will, however, not contradict Theorem \ref{maintheorem} below because $(M,g)$ is null (and timelike) geodesically incomplete.

On the other hand, for $k\geq 2$ and $n\geq 4$, any other round sphere in $\Sigma$ whose interior includes some of the $z_i$'s but excludes others will still be enclosing in the sense of definition \ref{asympticallyregular}. In this case $\overline{\Sigma}_+$ is not topologically $S\times [0,+\infty)$ as demanded in the original Gannon-Lee context.
\end{exe}

\begin{exe}
\label{morexe2}
Let $(M_0,g_0)$ be any smooth Riemannian $n$-manifold. On $M_0$, pick a smooth, real-valued, strictly positive function $\beta_0$, and a smooth 1-form $\omega_0 \in \Omega^1(M_0)$. Fix also a strictly positive smooth function $\Lambda_0 \in C^{\infty}(\mathbb{R}\times M_0)$. Then, the {\em standard conformastationary spacetime} associated with the data $(M_0,g_0,\beta_0,\omega_0,\Lambda_0)$ is $(M, g)$, where $M:= \mathbb{R} \times M_0$, and
\begin{equation}
\label{standardstationary}
g = \Lambda_0^2(-\beta ^2 d\pi_1 \otimes d\pi_1 + \omega \otimes d\pi_1 + d\pi_1 \otimes \omega + \pi_2^{\ast}g_0),
\end{equation}
where $\beta := \beta_0 \circ \pi_2$, $\omega := \pi_2^{\ast} \omega_{0}$, and $\pi_1$ [resp. $\pi_2$] is the projection of $M$ onto the $\mathbb{R}$ [resp. $M_0$] factor. The time-orientation of $(M,g)$ is chosen such that $\partial_t$, the lift to $M$ of the standard vector field $d/dt$ on $\mathbb{R}$, is future-directed. The vector field $\partial_t$ is then a timelike conformal Killing vector field. If $\Lambda_0\equiv 1$, then $(M,g)$ is said to be {\em standard stationary} (for the respective data), and if in addition $\omega_{0} \equiv 0$, then $(M,g)$ thus defined is said to be {\em standard static}.

The following facts about the standard conformastationary metric (\ref{standardstationary}) are germane for us here:
\begin{itemize}
\item[1)] The timelike conformal Killing vector field $\partial_t$ is complete, and hence Prop. \ref{criterion} applies and ensures that not only a standard conformastationary spacetime is reflecting, but all its coverings.
\item[2)] $\pi_1$ is a smooth {\em temporal function}, i.e., it has timelike gradient, and hence the hypersurfaces $\{t\} \times M_0$ are partial Cauchy hypersurfaces for each $t \in \mathbb{R}$. In particular, the spacetime is stably causal. By the general discussion in Section \ref{s1} they are causally continuous.
\item[3)] Clearly, the spacetime in Example \ref{morexe2} is standard static. Hence, in general conformastationary spacetimes are {\em not} causally simple.
\item[4)] $\partial _t$ is a piercing for, say, the partial Cauchy hypersurface $\Sigma:= \{0\} \times M_0$. Of course, in general neither $\Sigma$ is asymptotically regular, nor does $(M,g)$ satisfy the null convergence condition. But clearly there is a vast amount of topologies of $M_0$ and falloff conditions on the data $(g_0,\beta_0,\omega_0,\Lambda_0)$ that will ensure these properties.
\end{itemize}
\end{exe}

We are finally ready to state our main result.

\begin{thm}[{\bf Gannon-Lee theorem - new version}]
\label{maintheorem}
Let $(M,g)$ be an $n$-dimensional (with $n \geq 3$) null geodesically complete spacetime, which satisfies the null energy condition (i.e., $Ric(v,v)\geq 0$, for any null vector $v\in TM$), and possesses an asymptotically regular hypersurface $\Sigma \subset M$ pierced by some timelike vector field $X\in \mathfrak{X}(M)$. Let an enclosing surface $S \subset \Sigma$ be given, and assume, in addition, that at least one of the following conditions holds:
\begin{itemize}
\item[i)] $(M,g)$ and each one of its covering spacetimes $(\tilde{M},\tilde{g})$ are past reflecting, or else
\item[ii)] $S$ is simply connected and both $(M,g)$ and its {\em universal covering} are past reflecting.
\end{itemize}
Then, the group homomorphism $i_{\#}: \pi_1(S) \rightarrow \pi_1(\Sigma)$ induced by the inclusion $i: S \hookrightarrow \Sigma$ is surjective. In particular, if $S$ is simply connected, then so is $\Sigma$.
\end{thm}

In its version in \cite{costa_e_silva_gannon}, clauses $(i)$ and $(ii)$ were simply replaced by the condition that $(M,g)$ be {\em causally simple}. They obviously hold if $(M,g)$ is globally hyperbolic and $\Sigma$ is a Cauchy hypersurface, as in the original Gannon-Lee theorem. However, the proof in \cite{costa_e_silva_gannon} implicitly took for granted that causal simplicity holds on a certain Lorentzian {\em covering} of $(M,g)$, which as discussed in Example \ref{exe1} need not be the case\footnote{Indeed, note that the spacetime in Example \ref{exe1} does possess a piercing (via $\partial _t$) for its natural, $t$-constant partial Cauchy hypersurfaces. In particular, this shows that the existence of a piercing does not obviate this issue.}. Now, on the one hand this is {\em always} the case if $\Sigma$ is a Cauchy hypersurface as in the original Gannon-Lee theorems. On the other hand, in the present version, {\em provided one imposes past reflectivity on covering spacetimes as well as on $(M,g)$, these spacetimes might conceivably even be non-chronological}.

\bigskip
\noindent{\em Proof of Thm. \ref{maintheorem}}. \\
The key to the proof is establishing the following \\
\smallskip
\noindent{\em Claim:} {\em the interior region $\overline{\Sigma}_{-}$ in $\Sigma$ is compact.} (cf. \cite[Prop.\ 4.1]{costa_e_silva_gannon}). \\
\smallskip
\noindent In order to prove this Claim we begin by considering the set $\mathcal{H}^+$ of all the points $p$ of $E^+(S)(:= J^+(S)\setminus I^+(S))$ such that either $p \in S$ or else $p$ can be reached from $S$ by a future-directed null geodesic $\eta:[0,b]\rightarrow M$ with $\eta(0)\in S$, $\eta(b)=p$ and $\eta'(0) =K_{-}(\eta(0))$.

By a standard limit curve argument, either $\mathcal{H}^+$ is compact or else there exists a future-directed null geodesic $S$-ray $\gamma:[0,+\infty)\rightarrow M$, with $\gamma'(0) = K_{-}(\gamma(0))$. But the latter alternative is impossible due to the null convergence condition\footnote{Note that in this argument we can weaken our convergence assumption to an {\em averaged} null convergence condition in the following form: $\int_0^{+\infty}Ric(\gamma'(t),\gamma'(t))dt \geq 0$ along any null geodesic $\gamma:[0,+\infty)\rightarrow M$, with $\gamma'(0) = K_{-}(\gamma(0))$.} and the fact that $k_{-}>0$, which together imply the appearance of a focal point to $S$ along $\gamma$ incompatible with its maximal status \cite[Prop.\ 10.43]{oneill}. We conclude that $\mathcal{H}^+$ has to be compact.

Consider the closed set $T:=\partial I^+(\Sigma _+)\setminus \Sigma_+$. If we can show that $T\subset \mathcal{H}^+$, then it follows that $T$ is compact; in that case, arguing exactly as in the proof of Claim 3 in \cite[Prop.\ 4.1]{costa_e_silva_gannon}) we conclude that $\rho_X(T)\equiv \overline{\Sigma}_{-}$, where $\rho_X:M\rightarrow \Sigma$ is the retract associated with the piercing $X$ as discussed therein, and our Claim follows.

Suppose, then, by way of contradiction, that there is some $q \in T\setminus \mathcal{H}^+$. As mentioned in the Introduction, our strategy here is to adapt the proof of \cite[Thm.\ 2.3]{minguzzi_new} for the present context (see especially Fig. 3 in that reference). Let $(q_k)\subset I^+(q)$ be a sequence of points such that $q_k\rightarrow q$. Thus, $(q_k)\subset I^+(\Sigma_+)$. The maximal integral curve $\alpha$ of $X$ through $q$ must intersect $\Sigma$; it cannot do so to the future of $q$, or else this would violate the acausality of $\Sigma$. Thus, it either intersects $\Sigma$ at $q$ itself or to the past thereof. In the first case, since $q\notin \overline{\Sigma}_+$ we would have to have $q\in \Sigma_{-}$. In the second case, since $q\notin I^+(\overline{\Sigma} _+)$, we must have $q\in I^+(\Sigma_{-})$. In any case, $(q_k)\subset I^+(\Sigma_{-})$. Therefore, for each $k \in \mathbb{N}$ we have $I^-(q_k)\cap \Sigma_{\pm}\neq \emptyset$, whence we conclude that $I^-(q_k)\cap S\neq \emptyset$ since $I^-(q_k)\cap \Sigma$ is connected\footnote{To see this, just project a continuous curve in $I^-(q_k)$ between two points of $I^-(q_k)\cap \Sigma$ onto $\Sigma$ using $\rho_X$}. In other words, $(q_k)\subset I^+(S)$.

Fix a background complete Riemannian metric $h$ on $M$ with associated distance function $d_h$, and let $\sigma_k:[0,+\infty)\rightarrow M$ be a sequence of future-directed, future-inextendible timelike, $h$-arc-length-parametrized curves starting at $S$ and such that $\sigma_k(t_k)=q_k$ for some $t_k\in (0,+\infty)$. By the compactness of $S$ and the Limit Curve Lemma, we can, up to passing to subsequences, assume that $\sigma_k(0)\rightarrow r\in S$, and that there exists a future-directed, future inextendible $C^0$ causal curve $\sigma:[0,+\infty)\rightarrow M$ with $\sigma(0)=r$ such that
\[
\sigma_k|_C \rightarrow \sigma|_C
\]
$d_h$-uniformly for each compact set $C\subset [0,+\infty)$.

Suppose the sequence $(t_k)$ is bounded. Then we can assume, again up to passing to a subsequence, that $t_k\rightarrow t_0 \Rightarrow q=\sigma(t_0)$. Now, $q\notin S$, so $t_0>0$. By the achronality of $\partial I^+(\Sigma_+)$, the causal curve segment $\sigma |_{[0,t_0]}$ can be reparametrized as a future-directed null geodesic segment $\eta:[0,b]\rightarrow M$ without focal points to $S$ before $q=\eta(b)$. In particular, the null vector $\eta'(0)$ is normal to $S$, and hence it is either parallel to $K_{+}(r)$ or to $K_{-}(r)$. In the latter case, however, we'd have $q \in \mathcal{H}^+$, a contradiction. Thus we can assume, affinely reparametrizing $\eta$ if needed, that $\eta'(0)=K_{+}(r)$. But the acausality of $\Sigma$ implies that $\eta(0,b]$ cannot intersect $\Sigma$, in which case, as we discussed above, $q\in I^+(\Sigma_{-})$, and indeed the maximal integral curve $\alpha$ through $q$ intersects $\Sigma_{-}$. The continuous curve $\rho_X\circ \eta:[0,b]\rightarrow \Sigma$ enters initially in $\Sigma _{+}$, but $\rho_X\circ \eta(b) \in \Sigma_{-}$, so there exists some $s_0\in (0,b)$ for which $\rho_X\circ \eta(s_0) \in \Sigma _{+}$. But then $\eta(s_0) \in I^+(\Sigma _{+})$, so that $q\in I^+(\Sigma _{+})\cap \partial I^+(\Sigma_{+})$, again a contradiction. We conclude that $(t_k)$ must be unbounded.

We can assume $t_k\rightarrow +\infty$. If $\sigma:[0,+\infty)\rightarrow M$ never left $\partial I^+(S)$, it could again be reparametrized as a (future-complete) null geodesic $S$-ray $\gamma:[0,+\infty)\rightarrow M$ initially parallel to $K_{-}(\gamma(0))$, a contradiction. Thus, for some $b\in (0,+\infty)$ $\sigma(b)\in I^+(S)$ and we can pick an open set $U\ni \sigma(b)$ such that $U\subset I^+(S)$, and also pick $p \in I^-(\sigma(b),U)$. Since $\sigma_k(b)\rightarrow \sigma(b)$, and eventually $t_k>b$, eventually $p\ll \sigma_k(b)\leq \sigma _k(t_k)=q_k$, and we conclude that $q\in \overline{I^+(p)}$. But then, past reflectivity implies that $p\in \overline{I^-(q)}$. Since $p \in I^+(S)\subset I^+(\overline{\Sigma}_{+}) (\equiv I^+(\Sigma_{+}) )$, we thus have $q \in I^+(\Sigma_{+})\cap \partial I^+(\Sigma_{+})$. This final contradiction thus establishes the Claim.

The rest of the proof of Thm. \ref{maintheorem} now proceeds exactly as the proof of \cite[Thm.\ 2.1]{costa_e_silva_gannon}, only with the caveat that if clause $(ii)$ in the statement holds, then we may use the universal covering instead of the more elaborate one therein. Therefore, we omit further details here.
\qcd

\section*{Acknowledgments} IPCS is partially supported by the project MTM2016-78807-C2-2-P (Spanish MINECO with FEDER funds). Both authors wish to express their deep gratitude for the warm hospitality we received at the Dept. of Geometry and Topology of the University of Malaga, Spain, where this work was initiated. They extend their special thanks to Jos\'{e} Luis Flores and Miguel S\'{a}nchez for helpful discussions, as well as to Roland Steinbauer and Benedict Schinnerl for their careful reading of the manuscript and for their keen and kind comments thereon.

\end{document}